\documentclass[a4]{article}

 \usepackage{a4wide}
 \usepackage{latexsym}
 \usepackage{amsbsy}
 \usepackage{amssymb}
 \usepackage{amsmath}
 \usepackage[latin1]{inputenc}%needed for bibliography

\usepackage{color}

 \newtheorem{thm}{\hspace*{-1.1mm}}[section]
 \newcommand{\bt}{\begin{thm} {\bf Theorem }}
 \newcommand{\et}{\end{thm}}
 \newcommand{\bprop}{\begin{thm} {\bf Proposition }}
 \newcommand{\bcor}{\begin{thm} {\bf Corollary }}
 \newcommand{\blem}{\begin{thm} {\bf Lemma }}
 \newcommand{\bex}{\begin{thm} {\bf Example }\rm}
 \newcommand{\bexs}{\begin{thm} {\bf Examples }\rm}
 \newcommand{\bd}{\begin{thm} {\bf Definition }}
 \newcommand{\brem}{\begin{thm} {\bf Remark }\rm}
 \newcommand{\brs}{\begin{thm} {\bf Remarks }\rm}
 
 \newcommand{\pr}{\noindent{\bf Proof. }}
 \newcommand{\ep}{\nolinebreak{\hspace*{\fill}$\Box$ \vspace*{0.25cm}}}

 \renewcommand{\part}{Deo}

 \def\dj{d\kern-0.4em\char"16\kern-0.1em}
 \def\Dj{\mbox{\raise0.3ex\hbox{-}\kern-0.4em D}}

 \newcommand{\beq}{\begin{equation} }
 \newcommand{\eeq}{\end{equation} }

 \newcommand{\bea}{\begin{eqnarray}}
 \newcommand{\eea}{\end{eqnarray}}

 \newcommand{\beas}{\begin{eqnarray*}}
 \newcommand{\eeas}{\end{eqnarray*}}

 \newcommand{\beqs}{\begin{equation*}}
 \newcommand{\eeqs}{\end{equation*}}

 \newcommand{\bi}{\begin{itemize}}
 \newcommand{\ei}{\end{itemize}}

 \newcommand{\ben}{\begin{enumerate}}
 \newcommand{\een}{\end{enumerate}}

 \newcommand{\ba}{\begin{array}}
 \newcommand{\ea}{\end{array}}

 \newcommand{\ds}{\displaystyle}

 \newcommand{\R}{\mathbb R}
 \newcommand{\N}{\mathbb N}
 \newcommand{\C}{\mathbb C}

 \newcommand{\I}{\mathbb I}

 \newcommand{\cT}{\ensuremath{{\cal T}}}

 \newcommand{\la}{\langle}
 \newcommand{\ra}{\rangle}

 \newcommand{\al}{\alpha}
 \newcommand{\be}{\beta}
  
 \newcommand{\G}{\Gamma}

\newcommand{\Iat}{{}_a{\rm I}_t}
\newcommand{\Iot}{{}_0{\rm I}_t}

\newcommand{\Itaut}{{}_\tau{\rm I}_t}

\newcommand{\Jat}{{}_a{\rm J}_t}

\newcommand{\Itb}{{}_t{\rm I}_b}
\newcommand{\Itaub}{{}_{\tau}{\rm I}_b}
\newcommand{\ItT}{{}_t{\rm I}_T}

\newcommand{\Dat}{{}_a{\rm D}_t}
\newcommand{\Derot}{{}_0{\rm D}_t}

\newcommand{\Dtaut}{{}_{\tau}{\rm D}_t}

\newcommand{\Dtb}{{}_t{\rm D}_b}
\newcommand{\DtT}{{}_t{\rm D}_T}

\newcommand{\CDat}{{}_a^C{\rm D}_t}
\newcommand{\CDerot}{{}_0^C{\rm D}_t}

\newcommand{\CDTvt}{{}_{T_v}^C{\rm D}_t}

\newcommand{\CDtb}{{}_t^C{\rm D}_b}

\newcommand{\T}{\rm T}

\begin{document}
%%%%%%%%%%
%%%%%%%%%%

 \title{Control theory for nonlinear fractional dispersive systems}

 \author{Maja Joli\'c
         \footnote{Faculty of Sciences, University of Novi Sad, Serbia,
         Electronic mail: maja.jolic@dmi.uns.ac.rs}\\
         Sanja Konjik
         \footnote{Faculty of Sciences, University of Novi Sad, Serbia,
         Electronic mail: sanja.konjik@dmi.uns.ac.rs}\\
         Darko Mitrovi\'c
         \footnote{Faculty of Mathematics, University of Vienna, Oskar Morgenstern Platz-1, 1090 Wien, Austria,
         Electronic mail: darko.mitrovic@univie.ac.at}\thanks{Permanent address of D.M. is University of Montenegro, Montenegro}
      }

 \date{}
 \maketitle

 \begin{abstract}
 \noindent  We consider a terminal control problem for processes governed by a nonlinear system of fractional ODEs. 
 In order to show existence of the control, we first consider the linear counterpart of the system and reprove a number 
 of classical theorems in the fractional setting (representation of the solution through the Gramian type matrix, 
 Kalman's principle, equivalence of the controllability and observability). We are then in the position to use a fixed point theorem 
 approach and various techniques from the fractional calculus theory to get the desired result.   

 \vskip5pt
 \noindent
 {\bf MSC (2010):}
 %\begin{classcode}
  Primary: 26A33; Secondary: 34H05, 49J15 \\
 %\end{classcode}
 \end{abstract} 

%%%%%%%%%%
%%%%%%%%%%
 
%%%%% 
 \section{Introduction}
%%%%%

The recent development of computing technologies has allowed applications of more and more complicated and thus more 
accurate models for describing different phenomena in nature and society. In particular, this includes methods and approaches 
from fractional calculus which 
provide a significantly wider choice of models enabling us to fine--tune the equations describing the relevant phenomena.
Several classical partial differential equations (PDEs) have been used for many years as tools for modeling systems 
with applications in physics, biology, chemistry, medicine and engineering. On the other hand, recent experiments indicate 
that there is a large class of complex systems with different kinetics which have a microscopically complex behavior, 
and their macroscopic dynamics can not be described by the classical models.
To this end, we stress the results in \cite{DuWHu} where it was shown that fractional derivatives describe memory effects 
in different models and thus they appear to be unavoidable when it comes to the fine tuning of mathematical recast of natural phenomena.
More examples can be found in 
\cite{BagleyTorvik} where it has been shown that fractional calculus models of viscoelastic materials are consistent with 
the physical principles that govern such materials. In mathematical modeling in biology several papers regarding the 
fractional-order differential models of biological systems with memory, such as dynamics of tumor-immune system and
dynamics of HIV infection of CD4 + T cells and fractional-order predator-pray models have been considered in 
\cite{Ahmed, KutyniokLS, PriyaS}. Additional applications concern the spread of contaminants in underground water, 
network traffic, charge transport in amorphous semiconductors, cell diffusion processes, the transmission of signals through 
strong magnetic fields such as those found within confined plasma, etc. As for the purely mathematical point of view, 
most of the classical equations were successfully considered in the fractional setting. We shall avoid extensive listing of 
references and direct the reader to surveys with collections of applications given in e.g. \cite{BaleanuDiethelm, Podlubny}.

In the current paper, we plan to further expand the fractional calculus framework to nonlinear control theory. 
The concrete motivation for research is the model of interaction between the species which tend to avoid crowding. 
Unlike the situation from the seminal paper \cite{BGHP} where the space-time distribution of such species was considered, 
here we have a simpler case in the frame of which we are interested in a global density of species (it depends only on time). 
This means that instead of the PDE system from \cite{BGHP} we are interested in the ODE of the form 
\bea
\CDerot ^{\al} y(t) &=& -Af(y(t))y(t) + Bu(t),\quad t\in[0,T], \label{eq:prob1} \\
 y(0) &=& y_0, \label{eq:init_state}
 \eea
 where $\alpha\in (0,1)$, $y:[0,T]\to\mathbb{R}^d$ is the state of the system (in this case it is concentration of species),
 $u:[0,T]\to\mathbb{R}^N$ is the control function, $f:\R^d\to(0,\infty)$ is a continuous positive function, 
 $A\in\R^{d\times d}$ is a real symmetric positive semidefinite matrix, 
 $B\in\R^{d\times N}$ is a real matrix, $N\leq d$, and $y_0$ is the given initial state.
 
Non-positive definiteness of the matrix $A$ implies the dispersion tendency, while the function $f$ determines the intensity 
of dispersion (large $f(y)$ implies faster dispersion).
The matrix $B$ limits our possibility to influence the population densities, while $u$ represents controlled "birth" 
(it is appearance of new plants here) or "death" rate. We did not include explicitly the source term (representing "births" or "deaths"; 
see e.g. \cite{PT}) since we assume that we have total control over "births" and "deaths". However, this is merely a technical point 
and it is not difficult to include the additional term in the considerations. We stress that a fundamental issue here is how to choose a source term 
(usually called the control) in modeling equations which would govern the system from the given initial density to the prescribed final density. 

Finally, we note that the fractional derivative models memory effect, which in this case could be influence of the previously 
planted cultures on the ground quality and thus on our possibilities to cultivate new plants with a desired speed rate.         

%Such results have clearly a wide range of possible applications, 
%which makes control theory one of the most attractive fields in mathematics.

Unlike most investigations in the field of control theory so far, we shall address nonlinear problems. For that purpose, we shall
introduce a substantially new approach. The methodology 
for linear problems (numerical as well as analytic) is under intensive development and the existing techniques were adapted to 
the linear fractional situation (e.g. \cite{Biccari, BiccariHS, MicuZ-06}). In order to work with the nonlinear situation, we shall reduce the
problem to a linear one and use fixed point theorems and a priori estimates by adapting the method proposed in \cite{DKMN}.  
 
 Our aim is to show that for this class of the systems, under the assumption that $A$ and $B$ satisfy Kalman rank condition, 
 one can obtain global controllability. More precisely, that for any given final state $y_T$ there exists a control function 
 $u\in L^{2}([0,T];\R^N)$ such that the solution of the system \eqref{eq:prob1}, \eqref{eq:init_state} satisfies 
 \beq \label{eq:final_state}
 y(T)= y_T.
 \eeq 
 
 The idea is to linearize the problem by replacing the unknown function $f(y(t))$ in \eqref{eq:prob1} 
 with a fixed function  $f(v(t))$, where $v\in C([0,T];\R^d)$, and then use the results for linear control problems
 together with the fixed point theorem.
 
 We note here that the linear fractional control problems have been considered before in different settings (see \cite{BettayebD, LZ, Matignon}). 
 In particular, in \cite{LZ}, the authors prove that a fractional control system cannot be kept in the zero state if we cease with the control. 
 We stress that we have a different situation here, since we prove that the system is so called partially controllable, 
 i.e., that we are able to bring it into the final state (but not to keep it there without additional control). 
 
 The paper is organized as follows. After the Introduction, we have the section Preliminaries where we introduce necessary 
 notions and notations from the fractional calculus and functional analysis. In Section \ref{sec:LP}, we introduce the 
 linearized variant of \eqref{eq:prob1} and derive some controllability and observability results for linear time-varying 
 fractional control problems. In the final section, we prove the main results of the paper.
 
%%%%% 
 \section{Preliminaries}
%%%%%

Let us introduce basic notions regarding fractional derivatives, fractional integrals and their properties 
that will be needed for further considerations. Since in our problem the order of fractional derivatives is between
$0$ and $1$, we shall focus on that case. In addition, we shall assume enough regularity for functions that
are under consideration, so that all objects and identities are well defined (e.g. absolute continuity or
continuous differentiability). 

First, we recall the definition of the Riemann-Liouville fractional integral.

\bd 
The left and right Riemann-Liouville fractional integrals of order $\al\in(0,1)$ are given by 
$$
\Iat^{\alpha} f(t)=\frac{1}{\Gamma(\alpha)}\int\limits_{a}^{t} (t-\tau)^{\alpha-1}f(\tau)\,d\tau
$$
and 
$$
\Itb^{\alpha} f(t)=\frac{1}{\Gamma(\alpha)}\int\limits_{t}^{b} (\tau-t)^{\alpha-1}f(\tau)\,d\tau,
$$
where $\Gamma$ denotes the Gamma function.
\et

The Riemann-Liouville fractional operators of differentiation are introduced as follows.

\bd 
The left and right Riemann-Liouville fractional derivatives of order $\al\in(0,1)$,  are given by 
$$
\Dat^{\alpha} f(t)=\ds\frac{d}{dt}\big(\Iat^{1-\alpha}f(t)\big)=
\ds\frac{1}{\Gamma(1-\alpha)}\frac{d}{dt} \int\limits_{a}^{t} \dfrac{f(\tau)}{(t-\tau)^{\alpha}}\,d\tau
$$
and 
$$
\Dtb^{\alpha} f(t)=-\frac{d}{dt}\big(\Itb^{1-\alpha}f(t)\big)=
\ds\frac{-1}{\Gamma(1-\alpha)}\frac{d}{dt} \int\limits_{t}^{b} \dfrac{f(\tau)}{(\tau-t)^{\alpha}}\,d\tau.
$$
\et

We shall also need the Caputo fractional derivative, hence we recall its definition.

\bd 
The left and right Caputo  fractional derivatives of order $\al\in(0,1)$,  are given by 
$$
\CDat^{\alpha} f(t)=\ds\Iat^{1-\al}\big(f'(t)\big)
=\ds\frac{1}{\Gamma(1-\alpha)} \int\limits_{a}^{t} \dfrac{f'(\tau)}{(t-\tau)^{\alpha}}\,d\tau
$$
and 
$$
\CDtb^{\alpha} f(t)=\ds\Itb^{1-\alpha}\big(-f'(t)\big)
=\ds\frac{-1}{\Gamma(1-\alpha)} \int\limits_{t}^{b} \dfrac{f'(\tau)}{(\tau-t)^{\alpha}}\,d\tau.
$$
\et

The following proposition gives the fractional integration by parts formula that relates the Riemann-Liouville and Caputo 
fractional derivatives.

\bprop \label{parint} 
For $\al\in(0,1)$,
$$
\int\limits_{a}^{b} f(t)(\CDat^{\al}g(t))\,dt = (\Itb^{1-\al}f(t)g(t))|_{a}^b + \int\limits_{a}^{b} \Dtb^{\al}f(t)g(t)\,dt.
$$
\et

 The following lemma is proved in \cite[Th.\ 1.4, 1.5 \& 1.6]{Ferreira-Sign}.

\blem \label{lem:sign}
Let $\al\in(0,1)$ and $a\in C([0,T];\R)$. 
\bi
\item[a)]  If $x$ is a solution of the initial value problem 
$\CDerot^{\al} x(t)= a(t)x$, $x(0)=x_0>0$, 
then $x(t)>0$. 
\item[b)]  If $x$ is a solution of the initial value problem 
$\Derot^{\al} x(t)= a(t)x$, $\Iot^{1-\al}x(t)|_{t=0}=x_0>0$,
then $x(t)>0$. 
\item[c)]  If $x$ is a solution of the initial value problem 
$\DtT^{\al} x(t)= a(t)x$, $\ItT^{1-\al}x(t)|_{t=T}=x_T>0$, 
then $x(t)>0$. 
\ei
\et
 
 Following the lines of the proof of \cite[Prop.\ 3]{GallegosDM}, one can derive the next comparison result.
 
 \blem \label{lem:comp} 
 Let $\al\in(0,1)$ and $a_1, a_2\in C([0,T];\R)$ such that  $a_1(t)\leq a_2(t)$, $t\in[0,T]$.
 \bi
 \item[(i)] If $x_i:[0,T]\to\R$, $i=1,2$, is a solution of the initial value problem 
 $\CDerot^{\al} x_i(t)= a_i(t)x_i$, $x_i(0)=x_0$, 
 then $|x_1(t)|\leq |x_2(t)|$.
 \item[(ii)] If $x_i:[0,T]\to\R$, $i=1,2$, is a solution of the initial value problem 
 $\Derot^{\al} x_i(t)= a_i(t)x_i$, $\Iot^{1-\al}x_i(t)|_{t=0}=x_0$,
 then $|x_1(t)|\leq |x_2(t)|$.
 \item[(iii)] If $x_i:[0,T]\to\R$, $i=1,2$, is a solution of the initial value problem 
 $\DtT^{\al} x_i(t)= a_i(t)x_i$, $\ItT^{1-\al}x_i(t)|_{t=T}=x_T$,
 then $|x_1(t)|\leq |x_2(t)|$.
 \ei
 \et
 
 Let us recall the notion of the H\"{o}lder continuity. 
 A bounded continuous function $f:(a,b)\to\R^d$ is said to be the H\" older continuous  of order $\al\in(0,1]$ if there exists a constant $c_f>0$ such that  
 $$
 |f(x)-f(y)|\leq c_f |x-y|^{\al}, \text{ for every } x,y\in(a,b).
 $$
 The space of the H\" older continuous functions of order $\al\in(0,1]$ is denoted by $H^{\al}((a,b);\R^d)$, and it is a Banach space, when equipped with norm 
 $$
 \|f\|_{H^{\al}} = \sup\limits_{x\in(a,b)}|f(x)| + \sup\limits_{x,y\in(a,b), x\neq y}\frac{|f(x)-f(y)|}{|x-y|^{\al}}.
 $$
 
 Using the properties of the fractional integration operator $\Iat^{\al}$ given in \cite[Thm.\ 3.6]{SKM}, 
 we were able to prove the following property, see \cite[Lemma 4]{JKM}:
 
 \blem \label{lem:equicont} 
 Let $x_0\in\R^d$ and let $X$ be a subset of $C([0,T];\R^d)$ with the following properties:
 \begin{itemize}
 	\item[(i)] for every $x\in X$, $x(0)=x_0$;
 	\item[(ii)] there exists a constant $C>0$ such that for every $x\in X$, $\sup\limits_{t\in(0,T)}|\CDerot^{\al}x(t)|\leq C$.
 \end{itemize}
 Then $X$ is bounded in $H^{\al}((0,T);\R^d)$ and an equicontinuous subset of $C([0,T];\R^d)$.
 \et

%%% 
 \subsection{Mittag-Leffler functions}
%%%
 
 The special functions, such as the Gamma function, the Beta function, the Mittag-Leffler functions, etc, are an essential part of fractional calculus. 
 Here, we recall definition of the Mittag-Leffler functions. For a detailed analysis of this class of functions we refer to \cite{MitagLeffler-knjiga}.
 
 \bd 
 For $\al\in\C$, the one-parameter Mittag-Leffler function $E_{\al}(z)$ is defined by the series
 $$
 E_{\al}(z)=\sum\limits_{k=0}^{\infty}\frac{z^k}{\G(\al k+1)},\quad z\in\C.
 $$
 For $\al\in\C$, $\rm{Re}\{\al\}>0$, and $\be\in\C$, the two-parameter Mittag-Leffler function $E_{\al,\be}(z)$ is defined by the series
 $$
 E_{\al,\be}(z)=\sum\limits_{k=0}^{\infty}\frac{z^k}{\G(\al k+\be)},\quad z\in\C.
 $$ 
 \et
 
 The one-parameter Mittag-Leffler function is a special case of the two-parameter Mittag-Leffler function with $\beta=1$, i.e., 
 $E_{\al,1}(z)= E_{\al}(z)$.

%%% 
 \subsection{Fractional linear systems}
%%%

In this section we consider linear systems of fractional differential equations with time-varying coefficients 
and use results from \cite{Matychyn} to give an analytical representation of the solution. 

We start with the system involving the Riemann-Liouville fractional derivative, given by
\bea
\Dat^{\al}x(t) &=& A(t)x(t) + p(t), \quad t\in (a,b), \label{lin-RL}\\
\Iat^{1-\al}x(t)|_{t=a} &=& x_0, \label{lin-RL-pu}
\eea
where $A:[a,b]\to\R^{d\times d}$ is a continuous matrix valued function, $p:[a,b]\to\R^d$ is a continuous vector function, 
and $\al\in(0,1)$.

The state-transition matrix of this system is defined by a generalized Peano-Baker series 
\beq \label{PBS-RL}
\Phi(a,t)=\sum\limits_{k=0}^{\infty}\Iat^{k\circ\al}A(t),
\eeq
where $\ds\Iat^{0\circ\al}A(t)=\frac{(t-a)^{\al-1}}{\Gamma(\al)} \I$ and 
$\ds\Iat^{k\circ\al}A(t)=\Iat^{\al}(A(t) \Iat^{(k-1)\circ\al}A(t))$, for $k\geq1$ ($\I$ denotes identity matrix).

If the series in \eqref{PBS-RL} converges uniformly on $[a,b]$, then the state-transition matrix $\Phi(a,t)$ is the solution to initial value problem 
\beq \label{matr_prob_RL}
\Dat^{\al}\Phi(a,t) = A(t)\Phi(a,t), \ \  \Iat^{1-\al}\Phi(a,t)|_{t=a}= \I,
\eeq
and the solution of the system \eqref{lin-RL}, \eqref{lin-RL-pu} is given by (cf.\ \cite[Th.\ 3]{Matychyn})
\beq \label{resRL}
 x(t) = \Phi(a,t)x_0 + \int\limits_{a}^{t} \Phi(\tau,t)p(\tau)\,d\tau.
\eeq

Two examples of state-transition matrices, important for our analysis, are the following:
\bi
\item If $A(t)=A$ is a constant matrix, then 
$$
\Phi(a,t)=(t-a)^{\al-1}E_{\al,\al}(A(t-a)^{\al})=
(t-a)^{\al-1}\sum\limits_{k=0}^{\infty}\frac{A^k(t-a)^{k\al}}{\Gamma(\al(k+1))}.
$$
\item If $A(t)=Ag(t)$, where $A$ is a constant matrix and $g:[a,b]\to\R$ is a continuous function, then $\Phi(a,t)=\sum\limits_{k=0}^{\infty}A^k\Iat^{k\circ\al}(g(t))$.
Using the same arguments as in the proof of \cite[\emph{Prop. 1}]{JKM}, we conclude that $\Phi(a,t)$ converges uniformly on $[a,b]$.
\ei

Now we move to the case of systems with the Caputo fractional derivative:
\bea
\CDat^{\al}x(t) &=& A(t)x(t) + p(t), \quad t\in (a,b), \label{lin-Cap}\\
x(a) &=& x_0, \label{lin-Cap-pu}
\eea
where $A:[a,b]\to\R^{d\times d}$ is a continuous matrix valued function, $p:[a,b]\to\R^d$ is a continuous vector function, and $\al\in(0,1)$.

Here, the definition of the state-transition matrix, $\Psi(a,t)$, given in a form of a generalized Peano-Baker series is the following:
\beq \label{PBS-Cap}
\Psi(a,t)=\sum\limits_{k=0}^{\infty}\Jat^{k\circ\al}A(t),
\eeq
with $\ds\Jat^{0\circ\al}A(t)= \I$ and 
$\ds\Jat^{k\circ\al}A(t)=\Iat^{\al}(A(t) \Jat^{(k-1)\circ\al}A(t))$, for $k\geq1$.

Under the assumption that the series in \eqref{PBS-Cap} converges uniformly on $[a,b]$, the matrix $\Psi(a,t)$ satisfies the initial value problem:
$$
\CDat^{\al}\Psi(a,t) = A(t)\Psi(a,t), \quad \Psi(a,a)= \I.
$$
Furthermore, the solution to \eqref{lin-Cap}, \eqref{lin-Cap-pu} is given by (cf.\ \cite[Th.\ 5]{Matychyn})
\beq \label{resCap}
x(t) = \Psi(a,t)x_0 + \int\limits_{a}^{t} \Phi(\tau,t)p(\tau)\,d\tau.
\eeq

Here are two important examples in the case of the Caputo derivative:
\bi
\item If $A(t)=A$ is a constant matrix, then $\ds\Psi(a,t)=E_{\al}(A(t-a)^{\al})=\sum\limits_{k=0}^{\infty}\frac{A^k(t-a)^{k\al}}{\Gamma(\al k+1)}$.
\item If $A(t)=Ag(t)$, then $\ds\Psi(a,t)=\sum\limits_{k=0}^{\infty}A^k\Jat^{k\circ\al}(g(t))$, and from \cite[\emph{Prop. 2}]{JKM} it converges uniformly on $[a,b]$.
\ei

As stated in \cite{JKM}, if $A(t_1)$ and $A(t_2)$ commute for every $t_1,t_2\in[a,b]$, then the matrix $\Phi(\tau,t)$ 
from the nonhomogeneous part of solutions \eqref{resRL} and \eqref{resCap} coincides with the state-transition matrix 
of the initial value problem with the right Riemann-Liouville fractional derivative:
\bea
\Dtaut^{\al}\Phi(\tau,t) &=& A(\tau)\Phi(\tau,t), \quad \tau\in(a,t) \label{PhidesniRL} \\
\Itaut^{1-\al}\Phi(\tau,t)|_{\tau=t} &=& \I, \label{PhidesniRL_init}
\eea
where the upper limit $t$ is a fixed value from the interval $(a,b]$ and $\tau$ is the independent variable.

Furthermore, from \cite[Prop.\ 4.2.]{Gomoyunov} we have the following property of $\Phi(\tau,t)$.

\bprop 
Let $\Omega=\{(\tau,t)\in[a,b]^2| \tau\leq t\}$. The function $F(\tau,t):=(t-\tau)^{1-\al}\Phi(\tau,t)$ is continuous on $\Omega$ 
and there exists $H_{F}\geq 0$ such that
$$
\| F(t_1,s_1)-F(t_2,s_2)\|\leq H_{F}(|t_1-t_2|^{\al}+|s_1-s_2|^{\al}).
$$
\et

For further analysis it will be more convenient to consider $\Phi(\tau,t)$ as the solution of the right Riemann-Liouville problem 
of type \eqref{PhidesniRL}, \eqref{PhidesniRL_init}.
 
 \brem \label{rem:diag} (\textbf{Diagonalization})
 Since in our system the matrix $A$ is, by assumption, real and symmetric, there exist a diagonal matrix 
 $$
 D={\ensuremath{\mathrm{diag}}}(\lambda_1,\ldots,\lambda_d)
 $$ 
 and an orthogonal matrix $U$ such that 
 $$
 A=UDU^{-1}=UDU^{\T}.
 $$ 
 Elements on the main diagonal of $D$ are the eigenvalues of $A$. Using this diagonalization of $A$, 
 we will be able to diagonalize the system and derive some essential estimates of the solution.
\et

From \cite[Prop. 1 \& Rem. 3]{JKM} we have the following properties. 

\bprop \label{prop:Psi_Phi_1} 
Let $A$ be a real, symmetric and positive semidefinite matrix with eigenvalues $\lambda_1,\dots,\lambda_d$,  
$g\in C([a,b];[0,\infty))$ and let $\lambda=\max\limits_{1\leq i\leq d}\lambda_i$ and $M\geq0$ be such that $g(t)\leq M$, $t\in[a,b]$.
\bi
\item[(i)] If $\Psi(a,t)$ is a solution of the initial value problem 
\beq \label{matr_Cap_psi-} 
\CDat^{\al}\Psi(a,t) = -Ag(t)\Psi(a,t), \quad t\in(a,b), \quad \Psi(a,a)= \I,
\eeq
then $E_{\al}(-\lambda M(b-a)^{\al})\leq\|\Psi(a,t)\|\leq 1$.
\item[(ii)] Let $t\in(a,b]$. If $\Phi(\tau,t)$ is a solution of the initial value problem 
\beq \label{matr_phi-}
\Dtaut^{\al}\Phi(\tau,t) = -Ag(\tau)\Phi(\tau,t), \quad \tau\in(a,t), \quad \Itaut^{1-\al}\Phi(\tau,t)|_{\tau=t}= \I, 
\eeq 
then $E_{\al,\al}(-\lambda M(t-a)^{\al})\leq\|(t-\tau)^{1-\al}\Phi(\tau,t)\|\leq 1$.
\ei
\et

\bprop \label{prop:ogr_za_Psi_Phi} 
Let $A$ be a real symmetric matrix with eigenvalues $\lambda_1,\dots,\lambda_d$, $g\in C([a,b];\R)$, 
and let  $\ds\lambda_{\max}=\max\limits_{1\leq i\leq d}|\lambda_i|$ and $M\geq0$ be such that $|g(t)|\leq M$, $t\in[a,b]$.
\bi
\item[(i)] If $\Psi(a,t)$ is a solution of the initial value problem 
\beq\label{eq:matr_Cap_psi} 
\CDat^{\al}\Psi(a,t) = Ag(t)\Psi(a,t), \quad t\in(a,b), \quad \Psi(a,a)= \I,
\eeq
then $E_{\al}(-\lambda_{\max}M(b-a)^{\al})\leq\|\Psi(a,t)\|\leq E_{\al}(\lambda_{\max} M(b-a)^{\al})$.
\item[(ii)] Let $t\in(a,b]$. If $\Phi(\tau,t)$ is a solution of the initial value problem 
\beq\label{eq:matr_phi}
\Dtaut^{\al}\Phi(\tau,t) = Ag(\tau)\Phi(\tau,t), \quad \tau\in(a,t), \quad \Itaut^{1-\al}\Phi(\tau,t)|_{\tau=t}= \I, 
\eeq 
then $E_{\al,\al}(-\lambda_{\max}M(t-a)^{\al})\leq\|(t-\tau)^{1-\al}\Phi(\tau,t)\|\leq E_{\al,\al}(\lambda_{\max}M(t-a)^{\al})$.
\ei
\et

%%%
\subsection{A nonlinear fractional system}
%%%

A nonlinear system of fractional differential equations, associated with our control problem, is analyzed in \cite{JKM}. 
Here we recall the main result (cf. \cite[Th. 2]{JKM}) which will be substantial for construction of the solution of the control problem.

\bt\label{thm:particular}
Let $A\in\R^{d\times d}$ be a real symmetric positive semidefinite matrix, and $f:\R^d\to[0,\infty)$ 
a continuous function. Then, the Cauchy problem 
\beq \label{particular}
\begin{split}
	\CDerot^{\al}z(t) =& -Af(z(t))\, z(t), \quad t\in [0,T], \\
	z(0) =&\ z_0, 
\end{split}
\eeq
 has a solution $z\in C([0,T];\R^d)$.\\
Moreover, if $f$ is such that $F(z)=f(z)z$ is continuously differentiable on $\R^d$, then $z\in C([0,T];\R^d)\cap C^1((0,T];\R^d)$ 
is a unique solution and  $z'(t)=O(t^{\al-1})$, as $t\to 0$.
\et

\brem \label{rem:z_0} 
From the proof of \cite[Th. 2]{JKM} we have that the solution satisfies $|z(t)|\leq |z_0|$.
\et

\brem 
The uniqueness follows from \cite[Lemma 6.19.]{Diethelm} and the property that $z$ satisfies \eqref{particular} if and only if it satisfies 
the Volterra integral equation
$$
z(t)=z_0 - \frac{A}{\G(\al)}\int\limits_0^t(t-s)^{\al-1}f(z(s))z(s)\,ds.
$$
\et

%%%%%
 \section{Linearized control problem}
 \label{sec:LP}
%%%%%
   
 Let $A\in\R^{d\times d}$ be a real symmetric matrix, $B\in\R^{d\times N}$ a real matrix,  $g\in C([a,b];\R)$ and assume that $g(t)\neq 0$, 
 for almost every $t\in[a,b]$. Let $M>0$ be a fixed constant such that $|g(t)|\leq M$, for every $t\in[a,b]$. 
 In this section we are going to analyze the following linear control problem:
 \bea 
 \CDat ^{\al}y(t) &=& Ag(t)y(t) + Bu(t),\quad t\in[a,b], \label{eq:prob-lin} \\
 y(a) &=& y_0, \label{eq:init_state2}
 \eea
 which in the case $g(t)=f(v(t))$ coincides with the linearization of the problem \eqref{eq:prob1}, \eqref{eq:init_state}.
 
From the previous section we have that the solution of \eqref{eq:prob-lin}, \eqref{eq:init_state2} is given by
 \beq \label{resenje}
 y(t) = \Psi(a,t)y_0 + \int_{a}^{t} \Phi(\tau,t)Bu(\tau)\,d\tau, 
\eeq 
where $\Phi(\tau,t)=\sum\limits_{k=0}^{\infty}A^k\Itaut^{k\circ\al}(g(\tau))$.

%%% 
 \subsection{Controllability}
%%%
 
 Controllability of finite-dimensional linear fractional differential systems with time-invariant coefficients has been studied 
 in \cite{BalachandranTrujillo, BalachandranK, BettayebD, Matignon}. The obtained results show that controllability is equivalent 
 to invertibility of the Gramian matrix as well as to the Kalman rank condition. In this section, 
 we will derive analogous results for the system \eqref{eq:prob-lin}, \eqref{eq:init_state2}.

 To begin with, we shall first introduce the notion of controllability.
 
 \bd
 System \eqref{eq:prob-lin}, \eqref{eq:init_state2} is controllable if for every $y_b\in\R^d$
 there exists a control function $u:[a,b]\to\R^N$ such that the solution of the system satisfies $y(b)=y_b$.
 \et 
 
 In order to obtain the equivalent condition for controllability of the system \eqref{eq:prob-lin}, we introduce 
 the controllability Gramian matrix:
 \beq \label{eq:Gram}
 W(a,b) =\int\limits_{a}^{b}(b-t)^{1-\al}\Phi(t,b)BB^{\T}\Phi(t,b)^{\T}\,dt.
 \eeq
 Here $\Phi(t,b)$ denotes the matrix appearing in the nonhomogeneous part of the solution \eqref{resenje}. 
 More precisely, $\Phi(t,b)$ is the solution of the matrix problem \eqref{PhidesniRL}, \eqref{PhidesniRL_init} on the interval $[a,b]$, 
 with independent variable $t$ and the matrix of the system $A(t)=Ag(t)$. By $(\cdot)^{\T}$ we denote the transpose operator.
 Since $A$ is symmetric, we have $A^{\T}=A$, which implies that $\Phi(t,b)$ is symmetric too, i.e., $\Phi(t,b)^{\T}=\Phi(t,b)$. 
 
 \brem \label{rem:Gramian} 
 Compared to the Gramian matrix for integer-derivative control problems, in \eqref{eq:Gram} one could notice that the term $(b-t)^{1-\al}$ 
 is added. It was done in order to assure the convergence of the integral \eqref{eq:Gram}. 
 Let us note that it does not change the property of positive semidefiniteness of $W(a,b)$, since we have
 \beas
 x^{\T}W(a,b)x &=& \int\limits_a^b (b-t)^{1-\al}x^{\T}\Phi(t,b)B(t)B(t)^{\T}\Phi(t,b)^{\T}x\,dt \\
 &=& \int\limits_a^b (b-t)^{1-\al} |x^{\T}\Phi(t,b)B(t)|^2\,dt\geq 0,
 \eeas
 for any $x\in\R^d$. Hence, as in the integer-derivative case, the regularity of the Gramian is equivalent to its positive definiteness.

 Let us notice that using Remark \ref{rem:diag}, and introducing notation $\tilde{B}=U^TB$ and $x(t)=U^{\T}y(t)$, 
 we can reduce \eqref{eq:prob-lin}, \eqref{eq:init_state2} to an equivalent system: 
 \beq\label{eq:prob-lin-diag}
 \begin{split}
 	\CDat ^{\alpha}x(t)=&\ Dg(t)x(t) + \tilde{B}u(t), \quad t\in [a,b] \\
 	x(a) =&\ U^{\T}y_0.
 \end{split}
 \eeq
 Then, associated matrices $\Phi(\tau,t)$ and $\Phi_D(\tau,t)$ satisfy relation $\Phi(\tau,t)=U\Phi_D(\tau,t)U^{\T}$. 
 Furthermore, $\Phi_D(t,b)={\ensuremath{\mathrm{diag}}}(q_1(t,b),\ldots,q_d(t,b))$ where $q_i(t,b)$ is solution to
 $$
 \Dtb^{\al}q_i(t,b) = \lambda_i g(t)q_i(t,b), \quad \Itb^{1-\al}q_i(t,b)|_{t=b}=1, 
 \quad i=1,\dots,d.
 $$
 Hence, the controllability Gramians of \eqref{eq:prob-lin}  and \eqref{eq:prob-lin-diag}, denoted by $W(a,b)$ and $W_{D}(a,b)$, respectively, satisfy 
 \beas
 W(a,b) &=& \int\limits_a^b(b-t)^{1-\alpha}\Phi(t,b)BB^{\T}\Phi(t,b)^{\T}\,dt\\
 &=& \int\limits_a^b(b-t)^{1-\alpha}U\Phi_D(t,b)U^{\T}BB^{\T}U\Phi_D(t,b)^{\T}U^T\,dt\\
 &=& UW_{D}(a,b)U^{\T}.
 \eeas
 \et

 The following theorem provides an equivalent condition for controllability.
 
 \bt \label{thm:contr_equiv}
 The system \eqref{eq:prob-lin}, \eqref{eq:init_state2} is controllable if and only if the controllability Gramian matrix 
 $W(a,b)$ is nonsingular.
 \et
 
 \pr
 ($\Leftarrow$) If $W(a,b)$ is nonsingular, then the control 
 \beq \label{eq:u_*}
 u^*(t)= (b-t)^{1-\al}B^{\T}\Phi(t,b)^{\T}W(a,b)^{-1}[y_b-\Psi(a,b)y_0]
 \eeq
 is well defined and it steers the solution of the system to desired value $y_b$.
 Indeed, inserting $u=u^*$ in \eqref{resenje} we obtain
 \beas
 y(b) &=& \Psi(a,b)y_0 + \int\limits_{a}^{b} \Phi(\tau,b)Bu^*(\tau)\,d\tau \\
 &=& \Psi(a,b)y_0 + \int_{a}^{b} (b-\tau)^{1-\al} \Phi(\tau,b)B  B^{\T}\Phi(\tau,b)^{\T}W(a,b)^{-1}[y_b-\Psi(a,b)y_0] \,d\tau \\
 &=&  \Psi(a,b)y_0 + W(a,b)W(a,b)^{-1}[y_b-\Psi(a,b)y_0] \\
 &=& y_b.
 \eeas
 
 ($\Rightarrow$) Suppose that the system is controllable and that $W(a,b)$ is singular.
 Then there exists a column vector $w\in\R^d\setminus\{0\}$ such that $w^{\T}W(a,b)w=0$, which yields
 \beq \label{eq:gram0}
 \int\limits_{a}^bw^{\T}(b-t)^{1-\al}\Phi(t,b)BB^{\T}\Phi(t,b)^{\T}w\,dt=0.
 \eeq
 If 
 $$\phi(t)=w^{\T}\Phi(t,b)B$$
 then \eqref{eq:gram0} implies that $\phi(t)=0$, almost everywhere on $[a,b]$.
 Since the system is controllable, there exists a control $u$ such that the solution of the system satisfies $y(b)=w$, when $y(a)=0$.
 Then we have
 $$
 w=y(b)=\int\limits_{a}^{b}\Phi(t,b)Bu(t)\,dt,
 $$
 and 
 $$
 \|w\|^2=w^{\T}w=\int\limits_{a}^{b}w^{\T}\Phi(t,b)Bu(t)\,dt = \int\limits_{a}^{b}\phi(t)u(t)\,dt=0.
 $$
 This leads to a contradiction with $w\neq0$.
 \ep

 In the classical case (i.e., with integer order derivatives), as well as in the case of fractional systems with constant coefficients, 
 one can prove equivalence between non-singularity of the Gramian and the Kalman rank condition. 
 Next theorem provides the same result in the time-varying fractional context.
 
 \bt\label{thm:kalman-rank}
 System \eqref{eq:prob-lin}, \eqref{eq:init_state2} is controllable if and only if the Kalman rank condition is satisfied, 
 i.e., ${\rm rank}[B|AB|A^2B|\cdots| A^{d-1}B]=d$.
 \et
 
 \pr
  ($\Rightarrow$) Suppose that the system is controllable and that ${\rm rank}[B|AB|A^2B|\cdots| A^{d-1}B]<d$.
 Then there exists a column vector $q\in\R^d\setminus\{0\}$ such that 
 \beq\label{kalmanq}
 qB=0,\ qAB=0,\ \dots\  qA^{d-1}B=0.
 \eeq
 Furthermore, controllability implies that for $y_b=0$ there exists a control function $u_0$ such that the solution of the system satisfies $y(b)=0$.
 Then, from the representation of the solution \eqref{resenje}, we have
 $$
 0=y(b)=\Psi(a,b)y_0+\int\limits_a^b\Phi(\tau,b)Bu_0(\tau)d\tau.
 $$
 Therefore,
 \beq\label{Psi0}
 -\Psi(a,b)y_0 = \int\limits_a^b\Phi(\tau,b)Bu_0(\tau)\,d\tau = \int\limits_a^b \sum\limits_{k=0}^{\infty}A^kB(\Itaub^{k\circ \al}g(\tau))u(\tau)\,d\tau.
 \eeq
 From the Cayley-Hamilton theorem, it follows that for every $m\geq d$, $A^m=\sum\limits_{i=0}^{d-1}a_{m,i} A^{i}$, for some constants $a_{m,i}\in\R$. 
 Then, by induction, from \eqref{kalmanq} we conclude that for every $k\in\N$, $qA^kB=0$, and by multiplying \eqref{Psi0} by $q$, we obtain 
 $$
 -q\Psi(a,b)y_0 = \int\limits_a^b \sum\limits_{k=0}^{\infty}qA^kB(\Itaub^{k\circ \al}g(\tau))u(\tau)\,d\tau = 0.
 $$
 Since $\Psi(a,b)$ is regular matrix and $y_0$ is arbitrary, it follows that $q=0$ -- a contradiction.
 \\
 ($\Leftarrow$)
 Suppose that ${\rm rank}[B|AB|A^2B|\cdots| A^{d-1}B]=d$ and that $W(a,b)$ is singular. Then there exists $w\in\R^d\setminus\{0\}$  
 such that for almost every $t\in[a,b]$: 
 \beq \label{eq:gram-sing}
 w^{\T}\Phi(t,b)B=0.
 \eeq
 Applying $\Itb^{1-\al}$ to \eqref{eq:gram-sing}, and using $\Itb^{1-\al}\Phi(t,b)|_{t=b}=\I$, we obtain $w^{\T}B=0$.
 Moreover, from $\Dtb^{\al}\Phi(t,b)=Ag(t)\Phi(t,b)$ and \eqref{eq:gram-sing} it follows that
 $$
 w^{\T}Ag(t)\Phi(t,b)B=0, \quad \text{ for almost every } t\in[a,b].
 $$
 Since $g(t)$ is continuous and $g\neq0$ almost everywhere, we obtain $w^{\T}A\Phi(t,b)B=0$, for almost every $t\in[a,b]$,
 and therefore $w^{\T}AB=0$.
 
 By taking derivatives $\Dtb^{\al}$ of \eqref{eq:gram-sing} subsequently $k$-times, and then consecutively using $g\neq 0$, almost everywhere on $[a,b]$, 
 we get $w^{\T}A^k\Phi(t,b)B=0$, which together with $\Itb^{1-\al}\Phi(t,b)|_{t=b}=\I$ implies that $w^{\T}A^kB=0$, for every $k=0,1,2,...,d-1$.
 This contradicts the Kalman rank condition.
 \ep

%%%
 \subsection{Adjoint system. Observability}
%%%
 
 Controllability has its dual notion - observability, which is related to the corresponding adjoint system.
 For fractional linear systems, observability was introduced in \cite{Matignon}, where it was shown to be 
 equivalent to the positive definiteness of the observability Gramian matrix. 
 On the other hand, observability can be defined via the observability inequality, which indicates the existence 
 of a minimum point of a quadratic functional, which is closely related to the control function. 
 This approach was presented in \cite{Zuazua-14, Zuazua} for linear, time-invariant systems with integer derivatives, 
 and in \cite{MitrovicNU} for linear, time-invariant, parameter-depending systems with fractional derivatives. 
 In this section, we are going to adapt those notions to the framework of the Caputo (and in adjoint systems the Riemann-Liouville) fractional derivative, 
 and get results for the linear, time-varying systems of the form \eqref{eq:prob-lin}, \eqref{eq:init_state2}.

 We start by introducing the homogeneous adjoint system of \eqref{eq:prob-lin}, \eqref{eq:init_state2}. 
 We use $A^{*}$ to denote the adjoint operator of $A:\R^d\to\R^d$. Since $A$ is a real symmetric matrix, $A^*=A^{\T}=A$.
 
 Multiplying the equation $\CDat ^{\al}y(t)= Ag(t)y(t)$ by $z\in L^2([a,b];\R^d)$ and integrating over the interval $[a,b]$ gives us
 \beq \label{eq:adj1}
 \int\limits_a^b \la \CDat ^{\al}y(t),z(t)\ra \,dt= \int\limits_a^b \la Ag(t)y(t),z(t)\ra \,dt=
 \int\limits_a^b \la y(t), Ag(t)z(t)\ra \,dt.
 \eeq
 On the other hand, by applying Proposition \ref{parint} we have
 $$
 \int\limits_a^b \la \CDat ^{\al}y(t),z(t)\ra \,dt = 
 \la y(b),\Itb^{1-\al}z(t)|_{t=b}\ra - \la y(a),\Itb^{1-\al}z(t)|_{t=a}\ra + \int\limits_a^b \la y(t),\Dtb ^{\al}z\ra \,dt.
 $$
 This together with \eqref{eq:adj1} implies
 $$
 \int\limits_a^b \la y(t),\Dtb ^{\al}z(t)-Ag(t)z(t)\ra \,dt= 
 \la y_0,z_0\ra - \la y_b,z_b\ra, 
 $$
 where $z_b:=\Itb^{1-\al}z(t)|_{t=b}$ and $z_0:=\Itb^{1-\al}z(t)|_{t=a}$.
 
 Therefore, the adjoint system is given by:
 \bea
 \Dtb ^{\al}z(t) &=& Ag(t)\,z(t),\quad t\in[a,b], \label{eq:prob-adj} \\
 \Itb^{1-\al}z(t)|_{t=b} &=& z_b, \label{eq:init_adj}
 \eea
 with the solution $z(t)=\Phi(t,b)z_b$, where $\Phi(t,b)=\sum\limits_{k=0}^{\infty}A^k\Itb^{k\circ\al}(g(t))$ and
 $\Itb^{1-\al}\Phi(t,b)|_{t=b}=\I$.
 From Proposition \ref{prop:ogr_za_Psi_Phi} $(ii)$ it follows that the solution $z$ satisfies the following inequality
 \beq \label{ogr_za_z}
 (b-t)^{1-\al}|z(t)|\leq r_{\al+}|z_b|,
 \eeq
 where $r_{\al+}=E_{\al,\al}(\lambda_{\max}M(b-a)^{\al})$.

 Now we introduce the notion of observability.
 
 \bd \label{def:obs}
 System \eqref{eq:prob-adj}, \eqref{eq:init_adj} is observable if there exists $C>0$ such that for every $z_b\in\R^d$,
  solution of \eqref{eq:prob-adj}, \eqref{eq:init_adj} satisfies
 \beq \label{eq:obs-zT}
 |z_b|^2\leq C\int\limits_a^b|(b-t)^{1-\al}B^{\T}z(t)|^2\,dt.
 \eeq
 \et 
 
 \bprop 
 Observability inequality \eqref{eq:obs-zT} is equivalent to the following unique continuation property
 \beq \label{eq:obs_zT=0}
 B^{\T}z(t)=0,\ \text{ for almost every }t\in[a,b] \quad \Rightarrow \quad z_b=0.
 \eeq
 \et
 
 \pr 
 The implication \eqref{eq:obs-zT} $\Rightarrow$ \eqref{eq:obs_zT=0} follows directly.
 To prove \eqref{eq:obs_zT=0} $\Rightarrow$ \eqref{eq:obs-zT}, we define a mapping $\|\cdot\|_{r}:\R^d\to[0,\infty)$ such that for $z_b\in\R^d$,  
 $$\|z_b\|_{r}=\Big(\int\limits_a^b|(b-t)^{1-\al}B^{\T}z(t)|^2\,dt \Big)^{1/2},$$
 with $z(t)$ being a solution of the adjoint system 
 \eqref{eq:prob-adj}, \eqref{eq:init_adj}.
 Therefore, $z(t)=\Phi(t,b)z_b$ and \eqref{ogr_za_z} imply that $\|\cdot\|_r$ is well defined, and that for every $\mu\in\R$, $\|\mu z_b\|_r=|\mu|\|z_b\|_r$. 
 Furthermore, Minkowski inequality for $L^2([a,b];\R^d)$ provides the triangle inequality for $\|\cdot\|_r$.
 Hence, $\|\cdot\|_{r}$ is a seminorm on $\R^d$, and the assumption \eqref{eq:obs_zT=0} 
 implies that it is a norm. From the equivalence of all norms on $\R^d$ it follows that there exists a constant $C$ such that \eqref{eq:obs-zT} holds. 
 \ep
 
 The next step is to show that controllability and observability conditions are equivalent.  
 By multiplying equation \eqref{eq:prob-lin} with $z$ and integrating over the interval $[a,b]$ we obtain
 $$  
 \int\limits_a^b\langle \CDat^{\al}y,z\rangle \,dt = 
 \int\limits_a^b\langle Ag(t) y,z\rangle \,dt + \int\limits_a^b\langle Bu,z\rangle \,dt=
 \int\limits_a^b\langle y, Ag(t) z\rangle \,dt + \int\limits_a^b\langle Bu,z\rangle \,dt.
 $$
 On the other hand, Proposition \ref{parint} implies  
 $$
 \int\limits_a^b\langle \CDat^{\al}y,z\rangle \,dt =  
 \langle y,\Itb^{1-\al}z\rangle |_a^b + \int\limits_a^b\langle y,\Dtb ^{\al}z(t)\rangle \,dt,
 $$
 and hence
 $$
 \int\limits_a^b\langle y, \Dtb ^{\al}z(t)  - Ag(t)z\rangle \,dt =  
 -\langle y,\Itb^{1-\al}z\rangle |_a^b + \int\limits_a^b\langle Bu,z\rangle \,dt.
 $$
 Thus, solutions of \eqref{eq:prob-lin}, \eqref{eq:init_state2}, and the corresponding adjoint system \eqref{eq:prob-adj}, \eqref{eq:init_adj} satisfy 
 \beq \label{eq:eulerLagranze}
 \la y_0,z_0\ra - \la y_b,z_b\ra + \int\limits_a^b\langle u,B^{\T}z\rangle \,dt=0.
 \eeq

\bt 
System \eqref{eq:prob-lin}, \eqref{eq:init_state2} is controllable if and only if the adjoint system \eqref{eq:prob-adj}, \eqref{eq:init_adj} is observable.
\et

\pr 
($\Leftarrow$) Suppose \eqref{eq:prob-adj}, \eqref{eq:init_adj} is observable and 
\eqref{eq:prob-lin}, \eqref{eq:init_state2} is not controllable. Then the Gramian matrix $W(a,b)$ is singular and there exists 
$w\in\R^d\setminus\{0\}$ such that $B^{\T}\Phi(t,b)^{\T}w=0$, for almost every $t\in(a,b)$. 

Let $\tilde{z}=\Phi(t,b)^{\T}w$ be a solution to the adjoint system with initial condition $z_b=w$.
Then we have $B^{\T}\tilde{z}=0$, for almost every $t\in(a,b)$ and $z_b=w\neq0$ -- a contradiction with the observability.

($\Rightarrow$) Suppose that \eqref{eq:prob-lin}, \eqref{eq:init_state2} is controllable and that for every $C>0$ there exists 
some $z_b$ such that \eqref{eq:obs-zT} does not hold. Then we can define a sequence of numbers $\frac{1}{C_k}\to0$ 
and vectors $z_{kb}$ such that, without loss of generality, $|z_{kb}|=1$ and $|z_{kb}|^2>C_k\int\limits_a^b|B^{\T}(b-t)^{1-\al}z_k(t)|^2\,dt$.

Since  $z_{kb}$ is bounded sequence it has a convergent subsequence, which we do not relabel.
Then for $z_b=\lim\limits_{k\to\infty}z_{kb}$ we have $|z_b|=1$, and the solution $z(t)$ of the adjoint system with initial 
condition $z_b$ can be obtained as a limit of the solutions $z_k(t)$.

Now for solutions $z_k$ we have $0\leq\int\limits_a^b|B^{\T}(b-t)^{1-\al}z_k(t)|^2\,dt<\frac{1}{C_k}|z_{kb}|^2$.
Letting $k\to\infty$, we conclude that $B^{\T}z(t)=0$, almost everywhere on $[a,b]$.

Controllability assumption and \eqref{eq:eulerLagranze} imply that for $y_0=0$ and for every $y_b\in\R^d$ there exists a control $u$ 
such that 
$$
\la y_b,z_b\ra = \int\limits_a^b\langle u,B^{\T}z\rangle \,dt.
$$
Since $B^{\T}z(t)=0$ almost everywhere, it implies that $\la y_b,z_b\ra =0$, for every $y_b$. Therefore, $z_b=0$, leading to 
a contradiction with $|z_b|=1$.
\ep

 In the sequel we shall consider observability as a minimization problem. 
 More precisely, we shall show that observability inequality implies the coercivity 
 of the suitable quadratic functional, and that by minimizing that functional we obtain the optimal control $u$.
  
 Let us define the functional $J:\R^d\to \R$ by
 \beq \label{minimization}
 J(z_b)=\frac{1}{2}\int\limits_a^b |B^{\T}(b-t)^{\frac{1-\al}{2}}z(t)|^2 \,dt - \la y_b,z_b\ra + \la y_0,z_0\ra,
 \eeq where $z$ is the solution of \eqref{eq:prob-adj}, \eqref{eq:init_adj}.
 
 \bt \label{thm:minimizacija} 
If system \eqref{eq:prob-lin}, \eqref{eq:init_state2} is controllable (or equivalently \eqref{eq:prob-adj}, \eqref{eq:init_adj} is observable), 
 then the functional $J$ has a minimum. Moreover, the control function $u$ which steers the solution of \eqref{eq:prob-lin}, \eqref{eq:init_state2} 
 to the state $y(b)=y_b$ is given by 
 \beq \label{eq:u_preko_z}
 u(t)=B^{\T}(b-t)^{1-\al}\hat{z}(t),
 \eeq
 where $\hat{z}(t)$ is the solution of the adjoint system with final state $\Itb^{1-\al}\hat{z}(t)|_{t=b}=\hat{z}_b$ being the minimum point of $J$.
 \et
 
 \pr 
 Since $J$ is continuous and convex, it suffices to show that it is coercive, i.e., that $\ds\lim\limits_{|z_b|\to\infty}J(z_b)=\infty$.
 From \eqref{eq:obs-zT} it follows that
 $$
 |z_b|^2\leq C\int\limits_a^b|B^{\T}(b-t)^{\frac{1-\al}{2}}z(t)|^2(b-t)^{1-\al}\,dt
 \leq C(b-a)^{1-\al}\int\limits_a^b|B^{\T}(b-t)^{\frac{1-\al}{2}}z(t)|^2\,dt.
 $$
 Hence,
 $$
 J(z_b)\geq \frac{|z_b|^2}{2C(b-a)^{1-\al}} - \la y_b,z_b\ra + \la y_0,z_0\ra
 \geq\frac{|z_b|^2}{2C(b-a)^{1-\al}} - |\la y_b,z_b\ra| + |\la y_0,z_0\ra|,
 $$
 and, applying the Cauchy-Schwartz inequality, 
 $$
 J(z_b)\geq\frac{|z_b|^2}{2C(b-a)^{1-\al}} - |y_b||z_b| + |\la y_0,z_0\ra|.
 $$
 Here the right hand side tends to $+\infty$ when $|z_b|\to\infty$, implying that $\ds\lim\limits_{|z_b|\to\infty}J(z_b)=\infty$.
 
 Let $\hat{z}_b$ be the point where $J$ reaches its minimum. Then for every $z_b\in\R^d$ we have
 $$
 \lim\limits_{h\to0} \frac{J(\hat{z}_b+hz_b)-J(\hat{z}_b)}{h}=0.
 $$
 Calculating the limit above (see Appendix) we obtain 
 \beq \label{eq:uslov_za_z^}
 \int\limits_a^b\langle B^{\T}(b-t)^{\frac{1-\al}{2}}\hat{z}(t),B^{\T}(b-t)^{\frac{1-\al}{2}}z(t)\rangle \,dt + \la y_0,z_0\ra - \la y_b,z_b\ra =0.
 \eeq
 This is exactly \eqref{eq:eulerLagranze} for $u(t)=B^{\T}(b-t)^{1-\al}\hat{z}(t)$.
 \ep
 
 Theorem \ref{thm:minimizacija} provides several estimates of the control function.
 
 \bprop \label{prop:ogr_za_u}
 The control function $u$ which steers the solution of \eqref{eq:prob-lin}, \eqref{eq:init_state2} to the state $y(b)=y_b$ satisfies
 \bi
 \item[(i)] $|u(t)|\leq r_{\al+}\|B^{\T}\||\hat{z}_b|$, $t\in[a,b]$;
 \item[(ii)] $\|u\|_{L^2([a,b];\R^N)}\leq\sqrt{C}|y_b-\Psi(a,b)y_0|$, where $C$ is the observability constant from Definition \ref{def:obs}.
 \ei
 \et
 
 \pr  
 $(i)$ Follows directly from \eqref{ogr_za_z} and \eqref{eq:u_preko_z}.\\
 $(ii)$ Theorem \ref{thm:minimizacija} and \eqref{eq:eulerLagranze} imply
 \beq \label{eq:u_L2}
 \|u\|^{2}_{L^2([a,b];\R^N)}=\int\limits_a^b|(b-t)^{1-\al}B^{\T}\hat{z}(t)|^2\,dt=\la y_b,\hat{z}_b\ra-\la y_0,\hat{z}_0\ra.
 \eeq
 Let us show that $\la y_0,\hat{z}_0\ra=\la \Psi(a,b)y_0,\hat{z}_b\ra$.
 Suppose that $w$ is a solution of the initial value problem $\CDat^{\al}w(t)=Ag(t)w(t)$, $w(a)=y_0$. 
 Then $w(t)=\Psi(a,t)y_0$ and since $\hat{z}$ is a solution of the adjoint system, for every $t\in[a,b]$ we have
 $$
 \la\CDat^{\al}w(t),\hat{z}(t)\ra = \la Ag(t)w(t),\hat{z}(t)\ra= \la w(t), Ag(t)\hat{z}(t)\ra= \la w(t),\Dtb^{\al}\hat{z}(t) \ra.
 $$
 Therefore, $\ds\int\limits_a^b\la\CDat^{\al}w(t),\hat{z}(t)\ra \,dt =\int\limits_a^b\la w(t),\Dtb^{\al}\hat{z}(t) \ra \,dt$ 
 and Proposition \ref{parint} implies 
 $$
 \la w(b),\hat{z}_b\ra-\la w(a),\hat{z}_0\ra=0, \quad \mbox{ i.e., } \quad \la\Psi(a,b)y_0,\hat{z}_b\ra-\la y_0,\hat{z}_0\ra=0.
 $$
 Now, \eqref{eq:u_L2} reduces to $\|u\|^{2}_{L^2([a,b];\R^N)}=\la y_b-\Psi(a,b)y_0,\hat{z}_b\ra\leq |y_b-\Psi(a,b)y_0||\hat{z}_b|$, 
 and from \eqref{eq:obs-zT} it follows that $\|u\|^{2}_{L^2([a,b];\R^N)}\leq |y_b-\Psi(a,b)y_0|\sqrt{C}\|u\|_{L^2([a,b];\R^N)}$. After dividing by 
 $\|u\|_{L^2([a,b];\R^N)}$, we obtain the desired estimate.
 \ep

 \brem \label{rem:z^prekoW_c}
 From \eqref{eq:uslov_za_z^} and the expression for the solution of the adjoint system, it follows that for every $z_b\in\R^d$
 $$
 \int\limits_a^b\la B^{\T}(b-t)^{\frac{1-\al}{2}}\Phi(t,b)\hat{z}_b,B^{\T}(b-t)^{\frac{1-\al}{2}}\Phi(t,b)z_b\ra \,dt + \la y_0,z_0\ra - \la y_b,z_b\ra =0.
 $$
 Using the property $\la y_0,z_0\ra=\la \Psi(a,b)y_0,z_b\ra$, proved in the previous lemma, we obtain
 $$
 \la \int\limits_a^b (b-t)^{1-\al}\Phi(t,b)^{\T}BB^{\T}\Phi(t,b)\hat{z}_b\,dt, z_b\ra + \la\Psi(a,b)y_0,z_b\ra - \la y_b,z_b\ra =0.
 $$
 Since $A$ is symmetric, the state transition matrix $\Phi(t,b)$ will also be symmetric, and we have that the integral in the equation above
 is the controllability Gramian matrix. Therefore, for every $z_b\in\R^d$,
 $$
 \la W(a,b)\hat{z}_b + \Psi(a,b)y_0 - y_b, z_b \ra = 0,
 $$
 implying that $W(a,b)\hat{z}_b+ \Psi(a,b)y_0 - y_b=0$. Hence, we obtain the value of the minimum point of $J(z_b)$:
 \beq\label{eq:rem:z^prekoW_c}
 \hat{z}_b= W(a,b)^{-1}(y_b - \Psi(a,b)y_0).
 \eeq
 \et

%%%%%
\section{The main result}
%%%%%

In this section we address the question of global controllability of \eqref{eq:prob1}, \eqref{eq:init_state}. 
More precisely, we will prove that under the following conditions:
\bi
\item[(a1)] $A\in\R^{d\times d}$ is a real symmetric positive semidefinite matrix, $B\in\R^{d\times N}$ is a real matrix,
\item[(a2)] $A$ and $B$ satisfy the Kalman rank condition,
\item[(a3)] $f:\R^d\to(0,\infty)$ is continuous and further $F(y)=f(y)y$ is continuously differentiable on $\R^d$,
\ei
system \eqref{eq:prob1}, \eqref{eq:init_state} is controllable.
To that end, assume that $y_0$ is given initial state and $y_T$ is given final state, that we want to reach.
Let us denote by $z(t)$ the unique solution of the nonlinear initial value problem
\beq \label{nonlinear_z}
\begin{split}
	\CDerot^{\al}z(t) =& -Af(z(t))\, z(t), \quad t\in [0,T], \\
	z(0) =&\ y_0, 
\end{split}
\eeq
From Theorem \ref{thm:particular} we have: $z\in C([0,T];\R^d)\cap C^1((0,T];\R^d)$, 
\beq \label{ogr-z}
|z(t)|\leq |y_0|,
\eeq
and there exists $K_z>0$ such that for every $t\in(0,T]$ it holds
\beq \label{ogr-Kz}
|t^{1-\al}z(t)|\leq K_z.
\eeq

Next, for every $v\in C([0,T];\R^d)$ define constants $M_v, K_v$ and $T_v$ in the following way: 
\beq \label{def-MKT_v-frac}
M_v:=\max\limits_{t\in[0,T]}|f(v(t))|, \quad K_v:=\Big(\max\{1,M_v\}\Big)^{\frac{1}{\al}}, 
\quad \text{ and } \quad T_v:=T-\frac{T}{K_v^{\al}}.
\eeq
Notice that such $K_v$ satisfies 
\beq \label{ogr-f(v)_K_v^al}
\frac{\max\limits_{t\in[0,T]}|f(v(t))|}{K_v^{\al}}\leq 1.
\eeq
Further, we construct the solution $y$ in the following way: First, we let the system \eqref{eq:prob1} to be "uncontrolled", 
i.e., let $u=0$, up to the time $T_v\in(0,T)$, and then consider linearized control problem 
\beq \label{Nonlin.FDE.lineariz}
\begin{split}
	\CDerot^{\al} y(t)=&\ -Af(v)y + Bu, \quad t\in [T_v,T] \\
	y(T_v) =&\ z(T_v), \quad y(T)=y_T,
\end{split}
\eeq
and find its solution on $[T_v,T]$.
More precisely, for every $v\in C([0,T];\R^d)$ we define functions 
\beq\label{solution_y_u-frac}
y(t)= \begin{cases} z(t), &0\leq t \leq T_v\\
	y_2(t), & T_v< t \leq T
\end{cases},
\quad\quad
u(t)= \begin{cases} 0, &0\leq t \leq T_v\\
	u_2(t), & T_v< t \leq T
\end{cases},
\eeq
where  $z(t)$ is the restriction on $[0,T_v]$ of the solution to initial value problem \eqref{nonlinear_z}, and $y_2(t)$, $u_2(t)$ 
are the solution of the linear control problem \eqref{Nonlin.FDE.lineariz}. Notice that, although in \eqref{Nonlin.FDE.lineariz} 
we look for the solution on the interval $(T_v,T]$, the derivative $\CDerot^{\al}y(t)$ depends on the values of the solution 
on the whole interval $[0,t]$. Hence, we need to take into account the past of the system, 
i.e., the memory which is contained in the fractional derivative. By noticing
\beq \label{derivative+h}
\CDerot^{\al}y(t)=\frac{1}{\Gamma(1-\al)}\int\limits_{0}^{T_v}\frac{z'(s)}{(t-s)^{\al}}ds 
+ \frac{1}{\Gamma(1-\al)}\int\limits_{T_v}^{t}\frac{y_2'(s)}{(t-s)^{\al}}ds=h(t)+\CDTvt^{\al}y_2(t),
\eeq 
we can transform equation \eqref{Nonlin.FDE.lineariz} into
\beq \label{lineariz-h}
\begin{split}
	\CDTvt^{\al}y(t) =&\ -Af(v(t)) y(t) + Bu(t) - h(t) ,\quad t\in(T_v,T],\\
	y(T_v) =&\ z(T_v), \quad y(T)=y_T. 
\end{split}
\eeq
Then, we divide \eqref{lineariz-h} into two problems and find solution $y_2$ in the form $y_2(t)=y_p(t)+y_c(t)$, where $y_p$ solves 
\beq \label{prob-lin_p}
\begin{split}
	\CDTvt^{\al}y_p(t) =&\ -Af(v(t)) y_p(t) - h(t) ,\quad t\in(T_v,T],\\
	y_p(T_v) =&\ 0, 
\end{split}
\eeq
and $y_c$ and $u_2$ are solutions to 
\beq \label{prob-lin_c}
\begin{split}
	\CDTvt^{\al}y_c(t) =&\ -Af(v(t)) y_c(t) + Bu_2(t),\quad t\in(T_v,T],\\
	y_c(T_v) =&\ z(T_v)=y_{c,0}, \quad y_c(T)=y_T-y_p(T)=y_{c,T}. 
\end{split}
\eeq
From \eqref{ogr-z} we have 
\beq \label{y_c0}
|y_{c,0}|\leq |y_0|.
\eeq 
Denote by $\Psi_v(\tau,t)$ and $\Phi_v(\tau,t)$ the state transition matrices associated with the system \eqref{prob-lin_c}.
According to \eqref{resCap}, the solution to \eqref{prob-lin_p} is given by 
\beq \label{y_p}
y_p(t)=-\int\limits_{T_v}^t\Phi_{v}(\tau,t)h(\tau)d\tau.
\eeq

Using \eqref{ogr-Kz}, we get that for every $t\in(T_v,T]$, the function $h(t)$ from \eqref{derivative+h} satisfies:
\beas
|h(t)| &=&  \Big|\frac{1}{\Gamma(1-\al)}\int\limits_{0}^{T_v}\frac{z'(s)}{(t-s)^{\al}}ds\Big|
\leq  \frac{1}{\Gamma(1-\al)}\int\limits_{0}^{T_v}\frac{|s^{1-\al}z'(s)|}{s^{1-\al}(t-s)^{\al}}ds\\
&\leq& \frac{K_z}{\Gamma(1-\al)}\int\limits_{0}^{T_v}s^{\al-1}(t-s)^{-\al}ds \leq \frac{K_z}{\Gamma(1-\al)}\int\limits_{0}^{t}s^{\al-1}(t-s)^{-\al}ds.
\eeas
Introducing the change of variables $\xi=\frac{s}{t}$ in the above  integral, we obtain
\beq \label{h-z}
|h(t)|\leq \frac{K_z}{\Gamma(1-\al)}\int\limits_{0}^{1}\xi^{\al-1}(1-\xi)^{-\al}d\xi = \frac{K_z}{\Gamma(1-\al)}B(\al,1-\al)=K_z\Gamma(\al).
\eeq
Now, \eqref{y_p}, \eqref{h-z} and Proposition \ref{prop:Psi_Phi_1} $(ii)$ imply
\beq \label{ogr_za_y_p}
|y_p(t)|\leq \int\limits_{T_v}^t(t-\tau)^{\al-1}|h(\tau)|d\tau \leq K_z\Gamma(\al) \int\limits_{T_v}^t(t-\tau)^{\al-1}d\tau\leq K_z\Gamma(\al)\frac{T^{\al}}{\al}.
\eeq
Hence, for the final state $y_{c,T}$ from \eqref{prob-lin_c} we have
\beq \label{y_cT}
|y_{c,T}|\leq |y_T| +  \frac{K_z\Gamma(\al)T^{\al}}{\al} =: C_T.
\eeq

For the solution $y_c(t)$, we use properties of linear control presented in the previous section. We define $y_c(t)$ as the solution 
corresponding to the control $u_2(t)$ given by \eqref{eq:u_preko_z}. Using \eqref{eq:rem:z^prekoW_c}, we derive
\beq \label{control_u2}
u_2(t)=(T-t)^{1-\al}B^{\T}\Phi_v(t,T)^{\T}W_{v}^{-1}(y_T-\Psi_v(T_v,T)z(T_v)),
\eeq
where $W_{v}$ is the controllability Gramian associated to control problem \eqref{prob-lin_c}, i.e.,
$$
W_v=W_{v}(T_v,T)=\int\limits_{T_v}^{T} (T-t)^{1-\al}\Phi_v(t,T)BB^{\T}\Phi_v(t,T)^{\T}\,dt.
$$
Further, the solution $y_2$ is given by
\beq \label{solution_y_2-frac}
y_2(t)=\Psi_v(T_v,t)z(T_v) + \int\limits_{T_v}^t\Phi_v(\tau,t)Bu_2(\tau)\,d\tau.
\eeq

Let us prove some auxiliary results.

\blem \label{lem:convergence_MvKvTv} 
Let $f\in C(\R^d;(0,\infty))$ and let $\{v_n\}_{n\in\N}$ be a sequence in $C([0,T];\R^d)$ which converges uniformly on $[0,T]$ 
to a function $v\in C([0,T];\R^d)$. Then $\{f\circ v_n\}_{n\in\N}$ converges uniformly on $[0,T]$ to $f\circ v$. 
Moreover, $\{M_{v_n}\}_{n\in\N}$, $\{K_{v_n}\}_{n\in\N}$ and $\{T_{v_n}\}_{n\in\N}$ converge to $M_v$, $K_v$ and $T_v$, respectively.
\et

\pr 
Since $\{v_n\}_{n\in\N}$ is convergent sequence, it follows that it is bounded too, and there exists $K>0$ such that 
$$
\max\limits_{t\in[0,T]}|v(t)|\leq K \quad \text{and} \quad 
\max\limits_{t\in[0,T]}|v_n(t)|\leq K, \text{ for every } n\in\N.
$$
Let $B_K=\{x\in\R^d: |x|\leq K\}$. Since $f$ is continuous on $\R^d$ and $B_K$ is compact set, 
we have that $f$ is uniformly continuous on $B_K$. For a given $\varepsilon>0$, let $\delta>0$ be such that 
$$
(\forall x,y\in B_K)(|x-y|<\delta \Rightarrow |f(x)-f(y)|<\varepsilon).
$$
By choosing $n_0$ such that for $n>n_0$, $\max\limits_{t\in[0,T]}|v_n(t)-v(t)|<\delta$, we get
$$
n>n_0 \Rightarrow \max\limits_{t\in[0,T]}|f(v_n(t))-f(v(t))|<\varepsilon.
$$
Since $\varepsilon$ was arbitrary, it follows that $\{f\circ v_n\}_{n\in\N}$ converges uniformly to $f\circ v$ on $[0,T]$. 
Hence, $M_{v_n}=\|f\circ v_n\|_{C}\rightarrow\|f\circ v\|_{C}=M_v$, $n\to\infty$, which further implies $K_{v_n}\to K_v$ 
and $T_{v_n}\to T_v$, $n\to\infty$.
\ep

\blem \label{lem:LinFDE.ogr.za_Pshi_W_u} 
Let $v\in C([0,T];\R^d)$. There exist constants $\lambda\geq 0$, $c_w>0$ and $C_u>0$, not depending on $v$, such that:
\bi
\item[(i)] $\|\Psi_v(T_v,t)\|\leq 1$, for every $t\in[T_v,T]$;
\item[(ii)] $E_{\al,\al}(-\lambda M_v(t-\tau)^{\al})\leq \|(t-\tau)^{1-\al}\Phi_v(\tau,t)\|\leq 1$, for every $T_v\leq \tau\leq t\leq T$;
\item[(iii)] $\ds\|W_{v}^{-1}\|\leq\frac{K_v^{\al}}{c_w}$;
\item[(iv)] the control function given by \eqref{control_u2} satisfies: $|u_2(t)|\leq C_uK_v^{\al}$, for every $t\in[T_v,T]$.
\ei
\et

\pr 
Let $\lambda:=\max\{\lambda_1,\dots,\lambda_d\}$. Properties $(i)$ and $(ii)$ follow from Proposition \ref{prop:Psi_Phi_1}.

$(iii)$ Since we know that $W_{v}$ is nonsingular, we can define $\|W_{v}^{-1}\|=\frac{1}{s}$, where $s=\min\limits_{x\in S^{d-1}}|W_{v}x|$. 
Using Remark \ref{rem:Gramian}, we get that for every $x\in S^{d-1}$ 
$$
|W_{\al,v}x|\geq |x^{\T}W_{\al,v}x| = |x^{\T}UW_{v,D}U^{\T}x|.
$$
Since $U$ is orthogonal matrix, we have $\{x^{\T}U: x\in S^{d-1}\}=S^{d-1}$. Therefore, 
\beas
\min\limits_{x\in S^{d-1}}|W_{v}x|&\geq& \min\limits_{x\in S^{d-1}}|x^{\T}W_{v,D}x| \\
&=& \min\limits_{x\in S^{d-1}}\Big| \int\limits_{T_v}^T (T-t)^{1-\al} x^{\T}\Phi_{D,v}(t,T)\tilde{B}\tilde{B}^{\T}\Phi_{D,v}(t,T)^{\T}x\,dt \Big|\\
&=& \min\limits_{x\in S^{d-1}} \int\limits_{T_v}^T (T-t)^{1-\al}|x^{\T}\Phi_{D,v}(t,T)\tilde{B}|^2\,dt,
\eeas
where $\Phi_{D,v}(t,T)=diag(q_1(t),\dots,q_d(t))$, (cf.  Remark \ref{rem:Gramian}), with element $q_i(t)$ being the solution of
$$
\DtT^{\al}q_i(t) = -\lambda_if(v(t))q_i(t), \quad t\in[T_v,T], \quad \ItT^{1-\al}q_i(t)|_{t=T}=1.
$$
Using the same methods as in the proof of \cite[Prop.\ 1\ (ii)]{JKM}, we derive
$$
E_{\al,\al}(-\lambda_i M_v(T-t)^{\al}) \leq (T-t)^{1-\al}q_i(t)\leq 1, \quad t\in[T_v,T].
$$
Furthermore, from \eqref{def-MKT_v-frac}, we have
$$
M_v(T-t)^{\al}\leq M_v(T-T_v)^{\al}=M_v\frac{T^{\al}}{K_v^{\al}}\leq T^{\al}.
$$
Hence, $e_{\al}:=E_{\al,\al}(-\lambda T^{\al})\leq E_{\al,\al}(-\lambda_i M_v(T-t)^{\al})$, 
and we obtain uniform boundedness (with respect to both $v$ and $t$) of $q_i(t)$: 
\beq \label{ogr-q_i}
0<e_{\al}\leq (T-t)^{1-\al}q_i(t)\leq 1.
\eeq
Now, we have 
$$
\min\limits_{x\in S^{d-1}}|W_{v}x| \geq \min\limits_{x\in S^{d-1}} \int\limits_{T_v}^T (T-t)^{\al-1}|x^{\T}(T-t)^{1-\al}\Phi_{D,v}(t,T)\tilde{B}|^2\,dt,
$$
and from \eqref{ogr-q_i}, it follows that there exists a constant $c_x\geq0$, $c_x=c(U,e_{\al},B)$ (independent on $v$), such that
$$
\min\limits_{x\in S^{d-1}} \int\limits_{T_v}^T (T-t)^{\al-1}|x^{\T}(T-t)^{1-\al}\Phi_{D,v}(t,T)\tilde{B}|^2\,dt\geq c_x\int\limits_{T_v}^T(T-t)^{\al-1}\,dt = \frac{c_xT^{\al}}{K_v^{\al}}.
$$ 
Furthermore, the constant $c_x$ is strictly greater than $0$ since the assumptions (a1), (a2), (a3) imply positive definiteness 
of both $W_{v}$ and $W_{v,D}$ (cf. Proposition \ref{thm:kalman-rank} and  Remark \ref{rem:Gramian}). 

Denoting  $c_w=c_xT$, we obtain
$$
\min\limits_{x\in S^{d-1}}|W_vx|\geq \frac{c_w}{K_v^{\al}}.
$$
Hence, $\|W_v^{-1}\|\leq\frac{K_v^{\al}}{c_w}$.

$(iv)$ From \eqref{control_u2}, properties $(i)-(iii)$ and \eqref{ogr-z} we get 
\beas
|u_2(t)|&\leq& \|B^{\T}\|\|(T-t)^{1-\al}\Phi_v(t,T)^{\T}\| \|W_{v}^{-1}\|(|y_T| + \|\Psi_v(T_v,T)\||z(T_v)|) \\
&\leq& \frac{\|B^{\T}\|K_v^{\al}(|y_T|+|y_0|)}{c_w} = C_uK_v^{\al},
\eeas
with $C_u=\frac{\|B^{\T}\|(|y_T|+|y_0|)}{c_w}$.
\ep

Now, we are able to prove our main result.

\bt 
Assume that (a1), (a2) and (a3) hold. Then, for any $T>0$ and $y_0,y_T\in\R^d$, there exists $u\in L^2([0,T];\R^N)$ 
such that the solution of \eqref{eq:prob1}, \eqref{eq:init_state} satisfies $y(T)=y_T$.
\et 

\pr 
Define the mapping $\cT:C([0,T];\R^d)\to C([0,T];\R^d)$ which every $v\in C([0,T];\R^d)$ maps to the solution $y(t)$ given by 
 \eqref{solution_y_u-frac}, and constructed as described above. For the second part of the solution, given by  
 \eqref{solution_y_2-frac}, we have that \eqref{ogr-z} and Lemma \ref{lem:LinFDE.ogr.za_Pshi_W_u} imply that for every $t\in[T_v,T]$,
\beas
|y_2(t)| &=& | \Psi_v(T_v,t)z(T_v) + \int\limits_{T_v}^{t}\Phi_v(\tau,t)Bu(\tau)\,d\tau|\\
&\leq& | y_0| + \int\limits_{T_v}^{t}(t-\tau)^{\al-1}\|(t-\tau)^{1-\al}\Phi_v(\tau,t)\||Bu(\tau)|\,d\tau\\
&\leq& |y_0| + \|B\|C_uK_v^{\al}\int\limits_{T_v}^{t}(t-\tau)^{\al-1}\,d\tau
\leq |y_0| + \|B\|C_uK_v^{\al}\frac{(T-T_v)^{\al}}{\al}\\
&=&|y_0| + \|B\|C_uK_v^{\al}\left(\frac{T}{K_v}\right)^{\al} = |y_0| + \|B\|C_uT^{\al}=:C_y.
\eeas
Since on $[0,T_v]$ we have $|y(t)|=|z(t)|\leq |y_0|< C_y$, we conclude that
\beq\label{ogr-y-frac}
|y(t)|\leq C_y, \quad t\in[0,T].
\eeq

Now, let us show that the mapping $\cT$ is compact. Let $V$ be a bounded set in $C([0,T];\R^d)$. 
Since $f$ is continuous, there exists $K>1$ such that,  for every $v\in V$, 
\beq\label{ogr-f(v),f(z)}
\max\limits_{t\in[0,T]}|f(v(t))|\leq K^{\al} \quad \text{and}\quad \max\limits_{|z|\leq y_0}f(z)\leq K^{\al}.
\eeq 
Then, \eqref{def-MKT_v-frac}  implies $K_v\leq K$.
For the compactness of $\cT$, it suffices to prove that $Y:=\cT(V)$ is relatively compact set. 
To that end, let $\{y_n\}_{n\in\N}$ be a sequence in $Y$. For every $n\in\N$, denote by $v_n$ the function from $V$ such that $\cT(v_n)=y_n$.  
From \eqref{ogr-y-frac} we have that $\{y_n\}_{n\in\N}$ is uniformly bounded. Furthermore, for every $n$, the solution $y_n$ is equal to $z$ on $[0,T_v]$, 
i.e., $y_n$ satisfies \eqref{nonlinear_z} on $[0,T_v]$. Hence, from \eqref{ogr-f(v),f(z)} and \eqref{ogr-z} we obtain 
\beq \label{ogr-Dy-[0,T_v]}
|\CDerot^{\al}y_n(t)| = |-Af(z)z|\leq \|A\|K^{\al}|y_0|, \quad t\in[0,T_v].
\eeq
On the interval $(T_v,T]$ we have that $y_n(t)$ satisfies
$$
\CDerot^{\al} y_n(t) = -Af(v_n)y_n + Bu_{2,n}.
$$
Then, using \eqref{ogr-f(v),f(z)}, \eqref{ogr-y-frac} and Lemma \ref{lem:LinFDE.ogr.za_Pshi_W_u} $(iv)$, we get
\beq \label{ogr-Dy-[T_v,T]}
|\CDerot^{\al}y_n(t)|\leq \|A\|K^{\al}C_y + \|B\|C_uK_v^{\al}\leq (\|A\|C_y + \|B\|C_u)K^{\al}=:C_{\al} , \quad t\in(T_v,T].
\eeq
Therefore, \eqref{ogr-Dy-[0,T_v]} and \eqref{ogr-Dy-[T_v,T]} imply that the sequence of derivatives $\{\CDerot^{\al}y_n(t)\}_{n\in\N}$ 
is bounded on $(0,T)$, independently of $v$. Note here that $\CDerot^{\al}y_n(t)$ may not be continuous at $T_v$. 
Nonetheless, $\sup\limits_{t\in(0,T)}|\CDerot^{\al}y_n(t)|\leq C_{\al}$, for every $n\in\N$, and conditions of Lemma \ref{lem:equicont} are satisfied. 
Hence, $\{y_n\}_{n\in\N}$ is uniformly bounded and equicontinuous sequence in $C([0,T];\R^d)$, and by the Arzela-Ascoli theorem it follows that  
$\{y_n\}_{n\in\N}$ has a convergent subsequence. This concludes the proof of compactness.

Uniform boundedness of the solution \eqref{ogr-y-frac} also implies that the set 
$$
\{v\in C([0,T];\R^d): v=\omega\cT(v), \ \omega\in[0,1]\}
$$
is bounded.
It remains to show that $\cT$ is continuous. Assume that $\{v_n\}_{n\in\N}$ is a sequence in $C([0,T];\R^d)$ which converges uniformly 
to a function $\overline{v}\in C([0,T];\R^d)$. Then $\{v_n\}_{n\in\N}$ is bounded in $C([0,T];\R^d)$, and compactness of $\cT$ implies that 
the sequence $\{y_n\}_{n\in\N}=\{\cT(v_n)\}_{n\in\N}$ has a convergent subsequence $\{y_{n_k}\}_{k\in\N}$. Let $y=\lim\limits_{k\to\infty}y_{n_k}$. 
Now, from the construction of the solutions $y_{n_k}$, assumption that $v_n\to \overline{v}$, $n\to\infty$, uniformly on $[0,T]$, 
and Lemma \ref{lem:convergence_MvKvTv}, it follows that $y$ coincides with the solution $\overline{y}$, obtained for $v=\overline{v}$. 
Hence, $\{y_n\}_{n\in\N}$ converges to $y=\cT(\overline{v})$.

Therefore, $\cT $ satisfies conditions of the Leray-Schauder fixed point theorem, and we have the existence of the fixed point $y^*=v^*$, 
which is the desired solution of nonlinear control problem \eqref{eq:prob1}.
\ep

%%%%%
\section*{Appendix}
%%%%%

Here we derive \eqref{eq:uslov_za_z^} for the minimum point $\hat{z}_b$ of the functional $J$ from Theorem \ref{thm:minimizacija}. 

If $\hat{z}_b$ is the value for which $J$ reaches its minimum, then for every $z_b\in\R^d$ holds
\beq\label{J-var} 
\lim\limits_{h\to0} \frac{J(\hat{z}_b+hz_b)-J(\hat{z}_b)}{h}=0.
\eeq

Let $\hat{z}(t)$ and $z(t)$ be the solutions of the adjoint problem \eqref{eq:prob-adj} with $\Itb^{1-\al}\hat{z}|_{t=b}=\hat{z}_b$ 
and $\Itb^{1-\al}z|_{t=b}=z_b$, respectively.
Since $\hat{z}(t)=\Phi(t,b)\hat{z}_b$ and $z(t)=\Phi(t,b)z_b$, it follows that the solution $z_h$ of the equation \eqref{eq:prob-adj} 
that satisfies condition $\Itb^{1-\al}z_h|_{t=b}=\hat{z}_b+hz_b$ is given by $z_h(t)=\hat{z}(t)+hz(t)$.
Hence, using \eqref{minimization} we have
\beas
J(\hat{z}_b+hz_b) &=& \frac{1}{2}\int\limits_a^b |B^{\T}(b-t)^{\frac{1-\al}{2}}(\hat{z}(t)+hz(t))|^2 \,dt - \la yb,\hat{z}_b+hz_b\ra + \la y_0,\hat{z}_0+hz_0\ra \\
 &=& \frac{1}{2}\int\limits_a^b \la B^{\T}(b-t)^{\frac{1-\al}{2}}(\hat{z}(t)+hz(t)),B^{\T}(b-t)^{\frac{1-\al}{2}}(\hat{z}(t)+hz(t))\ra \,dt \\
 && - \la y_b,\hat{z}_b\ra - h\la y_b,z_b\ra + \la y_0,\hat{z}_0\ra + h\la y_0,z_0\ra \\
 &=& \frac{1}{2}\int\limits_a^b |B^{\T}(b-t)^{\frac{1-\al}{2}}\hat{z}(t)|^2 \,dt + h\int\limits_a^b \la B^{\T}(b-t)^{\frac{1-\al}{2}}\hat{z}(t),B^{\T}(b-t)^{\frac{1-\al}{2}}z(t)\ra \,dt\\
 && +  \frac{h^2}{2}\int\limits_a^b |B^{\T}(b-t)^{\frac{1-\al}{2}}z(t)|^2 \,dt - \la y_b,\hat{z}_b\ra - h\la y_b,z_b\ra + \la y_0,\hat{z}_0\ra + h\la y_0,z_0\ra \\
 &=& J(\hat{z}_b) + h\int\limits_a^b \la B^{\T}(b-t)^{\frac{1-\al}{2}}\hat{z}(t),B^{\T}(b-t)^{\frac{1-\al}{2}}z(t)\ra \,dt  - h\la y_b,z_b\ra + h\la y_0,z_0\ra\\
 && +  \frac{h^2}{2}\int\limits_a^b |B^{\T}(b-t)^{\frac{1-\al}{2}}z(t)|^2 \,dt
\eeas 
Therefore
\beas
\frac{J(\hat{z}_b+hz_b) - J(\hat{z}_b)}{h} &=&
\int\limits_a^b \la B^{\T}(b-t)^{\frac{1-\al}{2}}\hat{z}(t),B^{\T}(b-t)^{\frac{1-\al}{2}}z(t)\ra \,dt - \la y_b,z_b\ra + \la y_0,z_0\ra \\
 && +  \frac{h}{2}\int\limits_a^b |B^{\T}(b-t)^{\frac{1-\al}{2}}z(t)|^2 \,dt
\eeas
Letting $h\to 0$ and using \eqref{J-var} we obtain 
$$ 
\int\limits_a^b\langle B^{\T}(b-t)^{\frac{1-\al}{2}}\hat{z}(t),B^{\T}(b-t)^{\frac{1-\al}{2}}z(t)\rangle \,dt - \la y_b,z_b\ra + \la y_0,z_0\ra=0, 
$$
which is precisely \eqref{eq:uslov_za_z^}.

%%%%%
 \section*{Declarations}
%%%%%

\subsection*{Funding}
The authors M. Joli\'c and S. Konjik acknowledge financial support of the Ministry of Education, Science and Technological Development 
of the Republic of Serbia (Grant No. 451-03-68/2022-14/200125) and bilateral project SK-SRB-21-0028. 
The work of D. Mitrovi\'c is supported by the project P 35508 of the Austrian Science Fund.
\subsection*{Conflict of interest}
The authors declare that they have no conflict of interest.

%%%%%%%%%%%%%%%%%%%%%%%%%%%%%%%%%%%%%%%%%%%%%%%%%%%%%%%%%%%%%%%%%%
%%%%%%%%%%%%%%%%%%%%%%%%%%%%%%%%%%%%%%%%%%%%%%%%%%%%%%%%%%%%%%%%%%
% \chapter*{Conclusion}
 %\label{ch:conclusion}
 %\addcontentsline{toc}{chapter}{Conclusion}
%%%%%%%%%%%%%%%%%%%%%%%%%%%%%%%%%%%%%%%%%%%%%%%%%%%%%%%%%%%%%%%%%%
%%%%%%%%%%%%%%%%%%%%%%%%%%%%%%%%%%%%%%%%%%%%%%%%%%%%%%%%%%%%%%%%%% 

%%%%%%%%%%%%%%%%%%%%%%%%%%%%%%%%%%%%%%%%%%%%%%%%%%%%%%%%%%%%%%%%%%
%%%%%%%%%%%%%%%%%%%%%%%%%%%%%%%%%%%%%%%%%%%%%%%%%%%%%%%%%%%%%%%%%%

%%%%%%%%%%%%%%%%%%%%%%%%%%%%%%%%%%%%%%%%%%%%%%%%%%%%%%%%%%%%%%%%%%
 \newpage \pagestyle{empty}
%%%%%%%%%%%%%%%%%%%%%%%%%%%%%%%%%%%%%%%%%%%%%%%%%%%%%%%%%%%%%%%%%%  

%%%%%%%%%%%%%%%%%%%%%%%%%%%%%%%%%%%%%%%%%%%%%%%%%%%%%%%%%%%%%%%%%%
%%%%%%%%%%%%%%%%%%%%%%%%%%%%%%%%%%%%%%%%%%%%%%%%%%%%%%%%%%%%%%%%%%
%%%%%%%%%%%%%%%%%%%%%%%%%%%%%%%%%%%%%%%%%%%%%%%%%%%%%%%%%%%%%%%%%%
%%%%%%%%%%%%%%%%%%%%%%%%%%%%%%%%%%%%%%%%%%%%%%%%%%%%%%%%%%%%%%%%%%

\begin{thebibliography}{99}
 
 
\bibitem{Ahmed}
{Ahmed, E., El--Sayed, A., El--Saka, H.}
\newblock Equilibrium points, stability and numerical solutions of fractional--order predator--prey and rabies models.
\newblock {\em J. Math. Anal. App.}, {\bf 325}:542--553, 2007.


\bibitem{BagleyTorvik}
{Bagley, R. L., Torvik, P. J.}
\newblock A theoretical basis for the application of fractional calculus to viscoelasticity.
\newblock {\em J.\ Rheology}, {\bf 27}:201--210, 1983.


\bibitem{BalachandranTrujillo}
{Balachandran, K., Govindaraj, V., Ortigueira, M.D., Rivero, M., Trujillo, J.J. }
\newblock Observability and controllability of fractional linear dynamical systems.
\newblock {\em IFAC Proceedings Volumes}, {\bf 46}(1):893--898, 2013.


\bibitem{BalachandranK}
{Balachandran, K., Kokila, J.}
\newblock On the controllability of fractional dynamical systems.
\newblock {\em Int. J. Appl. Math. Comput. Sci.}, {\bf 22}(3):523--531, 2012.


\bibitem{BaleanuDiethelm}
{Baleanu, D., Diethelm, K., Scalas, E., Trujillo, J.J.}
\newblock {\em Fractional calculus: Models and numerical methods}, volume 3 of {\em
	Series on Complexity, Nonlinearity and Chaos}.
\newblock World Scientific, Hackensack, 2012.


\bibitem{BGHP} {M.~Bertsch, M.E.~Gurtin, D.~Hilhorst, L.A.~Peletier,}
\newblock  On interacting populations that disperse to avoid crowding: preservation and segregation.
\newblock {\em J.Math.Biology} {\bf 23} (1985), 1--13.

\bibitem{BettayebD}
{Bettayeb, M., Djennoune, S.}
\newblock New results on the controllability and observability of fractional dynamical systems.
\newblock {\em J. Vibr. Control}, {\bf 14}(9-10):1531--1541, 2008.


\bibitem{Biccari}
{Biccari, U.}
\newblock Internal control for a non-local Schr\" odinger equation involving the fractional Laplace operator.
\newblock {\em Evol. Equ. Control Theory}, doi: 10.3934/eect.2021014, 2021.


\bibitem{BiccariHS}
{Biccari, U., Hern\'andez-Santamar\'ia, V.}
\newblock Controllability of a one-dimensional fractional heat equation: theoretical and numerical aspects.
\newblock {\em IMA J. Math. Control}, {\bf 36}(4):1199--1235, 2018.


\bibitem{Diethelm}
{Diethelm, K.}
\newblock {\em The Analysis of Fractional Differential Equations}.
\newblock Springer, Heidelberg, 2004.


\bibitem{DKMN}
{Djordjevi\'c, J., Konjik, S., Mitrovi\'c, D., Novak, A.}
\newblock Global controllability for quasilinear nonnegative definite system of ODEs and SDEs.
\newblock {\em J. Optim. Theory Appl.}, {\bf 190}(1):316--338, 2021.


\bibitem{DuWHu}
{Du, M., Wang, Z., Hu, H.} 
\newblock Measuring memory with the order of fractional derivative. 
\newblock {\em Sci. Rep.} {\bf 3}(3431), 2013.


\bibitem{Ferreira-Sign}
{Ferreira, R.}
\newblock Sign of the solutions of linear fractional differential equations and some applications.
\newblock {\em Vietnam J. Math.}, https://doi.org/10.1007, 2021.


\bibitem{GallegosDM}
{Gallegos, J., Duarte-Mermoud, M.}
\newblock Boundedness and convergence on fractional order systems.
\newblock {\em J. Comput. Appl. Math.}, {\bf 296}:815--826, 2016.

\bibitem{Gomoyunov}
{Gomoyunov, M.}
\newblock On representation formulas for solutions of linear differential equations with Caputo fractional derivatives.
\newblock {\em Fract. Calc. Appl. Anal.}, {\bf 23}:1141--1160, 2020.


\bibitem{MitagLeffler-knjiga}
{Gorenflo, R., Kilbas, A., Mainardi, F., Rogosin, S.}
\newblock {\em Mittag-Leffler Functions, Related Topics and Applications}.
\newblock Springer, Heidelberg, 2014.


\bibitem{JKM}
{Joli\'{c}, M., Konjik, S., Mitrovi\'{c}, D.}
\newblock On solvability for a class of nonlinear systems of differential equations with the Caputo fractional derivative.
\newblock {\em Fract. Calc. Appl. Anal.}, {\bf 25}:2126?2138, 2022.


\bibitem{KutyniokLS}
{Kutyniok, G., Lim, W. Q., Steidl, G.}
\newblock Shearlets: theory and applications.
\newblock {\em GAMM-Mitt.}, {\bf 37}(2):259--280, 2014.

\bibitem{LZ} L\"u, Q., Zuazua, E.
\newblock On the lack of controllability of fractional in time ODE
and PDE.
\newblock {\em Math. Control Signals Syst.} (2016) 28:10


\bibitem{Matignon}
{Matignon, D., d'Andr\'{e}a-Novel, B.}
\newblock Some results on controllability and observability of finite-dimensional fractional differential systems.
\newblock {\em IMACS, IEEE-SMC Proceedings Conference, Lille, France}, 952--956, 1996.


\bibitem{Matychyn}
{Matychyn, I.}
\newblock Analytical solution of linear fractional systems with variable coefficients involving Riemann-Liouville and Caputo derivatives.
\newblock {\em Symmetry}, {\bf 11}(11):1366, 2019.


\bibitem{MicuZ-06}
{Micu, S., Zuazua, E.}
\newblock On the controllability of the fractional order parabolic equation.
\newblock {\em SIAM J. Control Optim.}, {\bf 44}(6):1950--1972, 2006.


\bibitem{MitrovicNU}
{Mitrovi\' c, D., Novak, A., Uzunovi\' c, T.,}
\newblock Averaged control for fractional ODEs and fractional diffusion equations.
\newblock {\em J. Funct. Spaces}, {\bf 2018}, Article ID 8095728, 2018.


\bibitem{Podlubny}
{Podlubny, I.}
\newblock {\em Fractional Differential Equations}, volume 198 of {\em
	Mathematics in Science and Engineering}.
\newblock Academic Press, San Diego, 1999.

\bibitem{PT} M.A.Pozio, A. Tesei, {\em Degenerate Parabolic Problems
in Population Dynamics}, Japan J. AppL Math. {\bf 2} (1985), 351--380.

\bibitem{PriyaS}
{Priya, G. Sudha, Prakash, P., Nieto, J.J., Kayar, Z.}
\newblock Higher-order numerical scheme for the fractional heat equation with Dirichlet and Neumann boundary conditions.
\newblock {\em Numer. Heat Transf. B: Fundam.}, {\bf 63}(6):540--559, 2013.


\bibitem{SKM} 
{Samko, S. G., Kilbas, A. A., Marichev, O. I.}
\newblock {\em Fractional integrals and derivatives}. 
\newblock  Gordon and Breach Science Publishers, Yverdon, 1993.


\bibitem{Zuazua-14}
{Zuazua, E.}
\newblock Averaged Control.
\newblock {\em Automatica}, {\bf 50}(12):3077--3087, 2014.


\bibitem{Zuazua}
{Zuazua, E.}
\newblock Controllability of partial differential equations. 
\newblock 3rd cycle. Castro Urdiales (Espagne), pp.311. cel-00392196, 2006.



\end{thebibliography}
\end{document}